\documentclass[12pt,psfig]{article}
\usepackage{amsmath, amsthm, amssymb}
\usepackage{graphicx,epsfig}
\usepackage{amsfonts}
\usepackage{amscd}
\usepackage{mathrsfs}
\usepackage{cite,calc,xcolor,subfigmat,mathrsfs,dsfont}
\usepackage{hyperref}
\usepackage{latexsym}

\newcommand{\ZS}{\mathscr{Z}}

\newcommand{\PM}{\mathscr{P}}

\newcommand{\E}{\mathscr{E}}

\newcommand{\RC}{\mathcal{R}}
\newcommand{\f}{\mathbb{F}}
\newcommand{\N}{\mathbb{N}}
\newcommand{\C}{\mathcal{C}}

\newcommand{\NC}{\mathscr{N}}
\newcommand{\RM}{\mathscr{R}}
\newcommand{\ta}{\mathscr{T}}

\newcommand{\NZ}{\mathbb{N}_0}

\newcommand{\Span}{{\rm span}}

\newcommand{\rank}{{\rm rank}}
\newcommand{\0}{\mathbf{0}}

\newcommand{\m}{\mathbf{m}}
\newcommand{\IM}{{\rm im}\,}
\newcommand{\ie}{{\em i.e.,} }
\newcommand{\eg}{{\em e.g.,} }

\newcommand{\ddx}{\frac{\partial}{\partial x}}
\newcommand{\ddy}{\frac{\partial}{\partial y}}
\newcommand{\ddz}{\frac{\partial}{\partial z}}
\newcommand{\ddY}{\frac{\partial}{\partial Y}}
\newcommand{\ddZ}{\frac{\partial}{\partial Z}}

\newtheorem{thm}{Theorem}[section]

\newtheorem{lem}[thm]{Lemma}

\theoremstyle{definition}
\newtheorem{defn}[thm]{Definition}
\newtheorem{exm}[thm]{Example}
\newtheorem{rem}[thm]{Remark}

\numberwithin{equation}{section}
\def\be {\begin{equation}}
\def\ee {\end{equation}}
\def\ba {\begin{eqnarray}}
\def\ea {\end{eqnarray}}
\def\bes {\begin{equation*}}
\def\ees {\end{equation*}}
\def\bas {\begin{eqnarray*}}
\def\eas {\end{eqnarray*}}
\def\bpr {\begin{proof}}
\def\epr {\end{proof}}

\begin{document}
\baselineskip=18pt
\renewcommand {\thefootnote}{\dag}
\renewcommand {\thefootnote}{\ddag}
\renewcommand {\thefootnote}{ }

\pagestyle{empty}

\begin{center}
                \leftline{}
                \vspace{-0.00 in}
{\Large \bf
Normal forms for Hopf-Zero singularities with nonconservative nonlinear part
} \\ [0.4in]

{\large Majid Gazor$^{*},$ Fahimeh Mokhtari}
\footnote{$^*\,$Corresponding author. Phone: (98-311) 3913634; Fax: (98-311) 3912602; Email: mgazor@cc.iut.ac.ir}


\vspace{0.15in}
{\small {\em Department of Mathematical Sciences,
Isfahan University of Technology
\\[-0.5ex]
Isfahan 84156-83111, Iran
}}

and
\vspace{0.15in}

{\large Jan A. Sanders}

\vspace{0.10in}
{\small {\em Department of Mathematics, Faculty of Sciences,
Vrije Universiteit, \\[-0.5ex]
Amsterdam 1081 HV, The Netherlands
}}

\noindent
\end{center}

\vspace{0.1in}
\baselineskip=15pt

\noindent \rule{6.5in}{0.012in}

\vspace{0.1in} \noindent

\begin{abstract}
In this paper we are concerned with the simplest normal form computation of the systems
\ba\label{0Eq1}
\dot{x}= 2x f(x, y^2+z^2),\, \dot{y}= z+y f(x, y^2+z^2),\, \dot{z}= -y+ z f(x, y^2+z^2),
\ea where \(f\) is a formal function with real coefficients and without any constant term. These are the classical normal forms of a larger family of systems with Hopf-Zero singularity. Indeed, these are defined such that this family would be a Lie subalgebra for the space of all classical normal form vector fields with Hopf-Zero singularity. The simplest normal forms and simplest orbital normal forms of this family with non-zero quadratic part are computed. We also obtain the simplest parametric normal form of any non-degenerate perturbation of this family within the Lie subalgebra. The symmetry group of the simplest normal forms are also discussed. This is a part of our results in decomposing the normal forms of Hopf-Zero singular systems into systems with a first integral and nonconservative systems.
\end{abstract}

\vspace{0.10in} \noindent {\it Keywords}: Parametric normal form; Hopf-Zero singularity; Eulerian vector fields.



\vspace{0.2in}

\section{Introduction}

Normal form theory is one of the most effective tools for the local study of nonlinear dynamical systems. The basic idea is to use permissible transformations and obtain a simplified vector field. Transformations are permissible that preserve certain dynamical features of the original system. The space of all permissible transformations form a group and acts on vector fields like an action of a group on a vector space. Consider a set of vector fields generated by the group acting on a given vector field. Then, the infinite level normal form of the vector field is to find a unique representative from this set. Thereby, the computation of infinite level normal forms is an important tool for classification of vector fields. The uniqueness of a normal form computation is determined by the (infinite) level of normal form through an specifically chosen normal form style and costyle. The level of a normal form assesses the remaining spectral data at our disposal in the permissible transformation space for further simplification of the system while a normal form style and costyle makes a unique choice for the normal form vector field in each level of normal form computation, see \eg \cite{baiderchurch,Gazor,GazorYuSpec,MurdBook,Sanders03}. Therefore, the infinite level normal form is sometimes called the simplest normal form or unique normal form when a normal form style have already been fixed. When a system has some symmetries, it is important that its normal form would preserve the symmetries. Although there are research results on the simplest normal forms of symmetric systems but they are considerably less than the existing results on normal forms without symmetry. There are several reasons for this. The first difficulty is to recognize the symmetries and then, to find the group of transformations preserving the symmetries. Therefore, one is also concerned with the space of the symmetric vector fields invariant under the group action. This needs a good understanding on the algebraic interactions of the symmetric and nonsymmetric vector fields with the transformation groups. Once all these are successfully accomplished, in most cases the normal form computation is more difficult in systems with symmetry than in systems without symmetry.

\pagestyle{myheadings}
\markright{{\footnotesize {\it M. Gazor, F. Mokhtari and J. A. Sanders
\hspace{2.3in} {\it Hopf-Zero singularity}}}}

Normal form decomposition of a nonsymmetric vector field into two symmetric vector fields can have many important potential applications.
For example Eulerian and Hamiltonian vector fields are two important families of vector fields. Therefore, the study of their dynamics is important. Normal form decomposition of arbitrary vector fields into Eulerian (nonconservative or dissipative) and Hamiltonian (conservative) vector fields are also important in both theory and applications. Wiggins \cite[Chapter 33]{Wiggins} remarks that transforming a system into an integrable Hamiltonian system plus a nonconservative perturbation facilitates ``a wealth of techniques for the global analysis of nonlinear dynamical systems such as Melnikov theory, perturbation theory for normally hyperbolic invariant manifolds, and Kolmogorov, Arnold, and Moser (KAM) theory''. Furthermore, Palaci\(\acute{\hbox{a}}\)n \cite{PalacChaos05} indicated that such kind of decomposition can be used in developing an ODE solver. For example using the obtained scalar function as the integral of a piece of the vector field in an ODE solver, can enhance its efficiency. These signify the importance of developing methods for such kind of decomposition.
Recently, some researchers have paid attention to this theory and have made important contributions to the subject, see \eg [3, 4, 7, 8, 18, 23--27, 29, 32].
Furthermore, the dynamics of a system may sometimes be analyzed by dynamics of its components. For example, an advection-diffusion model consists of advection terms and diffusion terms and its dynamics is well understood as a combination/competition of the dynamics associated with the advection and diffusion terms, see \eg \cite[Section 5.5]{Logan}. Therefore, the study of each components of a decomposed vector field may help to better understand the dynamics of the full system as a combined dynamics or as a competing behavior between its components. Thus, it is important to individually deal with the cases that the vector field is a quasi-Eulerian, or is a conservative vector field, see \cite{GazorMoazeni,GazorMokhtari1st}. In this paper a vector field is called {\em conservative} when it has a first integral and is called {\em nonconservative} when it does not have any first integral.

This paper deals with the nonconservative family of our upcoming results on such decomposition for Hopf-Zero singular vector fields.  Systems with Hopf-Zero singularity are important in applications.
There are several important research results on the simplest normal forms of Hopf-Zero singularity, see \cite{AlgabaHopfZ,ChenHopfZ03,ChenHopfZ,YuHopfZero}. However, there does not seem to exist any result on such kind of normal form decomposition for Hopf-Zero singularity. Although the Hamiltonian vector fields (symplectic structures) require an even dimensionality, the decomposition idea still can work. Indeed we instead work with conservative and nonconservative family of vector fields, see \cite{GazorMokhtari}. In this paper we are concerned with the study of nonconservative family. Indeed, we consider the family of Hopf-Zero singularities given by
\ba\label{Eq1}
v= z\ddy-y\ddz+ v_{[2]},
\ea where 
\ba
v_{[2]}:=f(x, y^2+z^2)E,\quad E:=x\ddx+\frac{1}{2}y\ddy+\frac{1}{2}z\ddz,
\ea and \(f\) is a formal power series function and without constant terms. Any result in this direction also helps in a better understanding of the more complicated algebraic structures involved in the three dimensional nilpotent singularities. In this paper the coefficients are taken from a field \(\f\) of characteristic zero, see \cite{Sanders03} where the coefficients are taken from a local ring containing the rational numbers. Since \(E\) is a quasi-Eulerian vector field, \(v_{[2]}\) is called a formal quasi-Eulerian vector field, that is a nonconservative vector field. The reader should note that this family of vector fields are indeed the classical (first level) normal form of a very much larger family of vector fields.
This, however, is independent from \cite{GazorMokhtari} and basically refers to an essentially separate problem. Our normal form computation is to use a transformation group that preserves the Eulerian structure and it would only simplify the function \(f.\) However, our results in \cite{GazorMokhtari} work with the bigger group, that is, the set of all near-identity transformations and may create conservative terms into the normal form system. Indeed, the quasi-Eulerian vector field \(E\) is defined such that its generated family would form a Lie algebra.

This paper is organized as follows. Section \ref{sec2} introduces a family of nonconservative vector fields that is a Lie algebra and presents a general notations for normal form computations. The simplest normal forms and orbital normal forms of this family of vector fields with non-zero quadratic part are obtained in Sections \ref{secNF} and \ref{secONF}, respectively. In Section \ref{secPNF} we prove that under certain technical conditions multi-parameter nonconservative perturbations of the vector field \(v\) given by Equation (\ref{Eq1}) can be transformed into
\ba\label{FinalPNF0}
w^{(\infty)}= z\ddy-y\ddz+ g(x, y^2+z^2, \mu)E, 
\ea where \(g(x,R,\mu)=x+\beta^0_1R+\sum^{r}_{i=0}R^{i}\mu_i+\beta^0_r R^{r}\) for a \(r\geq2.\)

\section{Nonconservative vector fields and notations}\label{sec2}

This section presents the necessary Lie algebra structure and notations that we need for normal form computations in the following sections.

\begin{defn} Define
\ba\label{Eul}
E^l_k&:=& x^{l+1}(y^2+z^2)^{k-l}\ddx+ \frac{1}{2} x^{l}
{(y^2+z^2)}^{k-l}y\ddy+\frac{1}{2} x^l{(y^2+z^2)}^{k-l}z\ddz,
(0\leq{l}\leq{k}).\quad
\ea Since \(E^l_k= x^{l}(y^2+z^2)^{k-l}E\) and \(E\) is a quasi-Eulerian vector field, we call \(E^l_k\) a quasi-Eulerian vector fields. Then define
\be\label{Lie}
\E:= \Span\Big\{\sum\beta^l_{k} E^l_k \;|\; \beta^l_{k}\in \f, {l+k\geq 1}, k\geq l\Big\}. \ee Let \(\Theta^0_0:= z\ddy-y\ddz.\) Therefore, the vector field given by Equation (\ref{Eq1}) can be represented by an element of \(\Theta^0_0+\E.\)
\end{defn} The following theorem follows from a straightforward computation.
\begin{thm}\label{iso} The vector space \(\{a\Theta^0_0+v|a\in \f, v\in\E\}\) is a Lie algebra with the following structure constants:
\bas\label{aa}
{[E^l_k, E^{m}_{n}]}&=& (n-k)E^{l+m}_{k+n},
\\
{[E_k^l, \Theta^0_0]}&=&0,
\eas for any \(l,k,m,n\in \NZ:=\N\cup\{0\}.\)
\end{thm}
\begin{proof} The space \(\E\) is a Lie algebra by the Lie bracket \([v,w]:= fE(g)E-gE(f)E\) for any \(v= f(x, y^2+z^2)E,w= g(x, y^2+z^2)E\in \E\) and \(E\) is considered as a differential operator acting on functions, see \cite[Page 2]{DumortierBook} and \cite{GazorMokhtari}. The structure constants follow from the fact that \(x^{l}(y^2+z^2)^{k-l}\) is an \(E\)-eigenfunction with eigenvalue \(k.\)
\end{proof} Note that one may consider the ring generated by polynomials \(x^{l}(y^2+z^2)^{k-l}.\) Then, it forms a Lie algebra; the Lie bracket of two eigenfunctions is given by their product multiplied with the substraction of their eigenvalues. Thus, normal form computations can be implemented on a computer algebra system using scalar valued functions rather than using vector fields. This enhances the efficiency of the computer program. The space \(\{a\Theta^0_0+v|a\in \f, v\in\E\}\) is indeed a Lie subalgebra of all the classical normal forms of Hopf-Zero singularity, see \eg \cite{AlgabaHopfZ,ChenHopfZ03,ChenHopfZ,MurdBook}.
\begin{thm}  There does not exist any first integral for any nonzero \(v\in \E,\) \ie any \(v\in \E\) is a non-conservative vector field.
\end{thm}

\bpr
Suppose that \(f=\sum_{m,n,k}\alpha_{m,n,k} x^my^nz^k\) is a first integral for \(gE\in \E\) where the function \(g\) is nonzero.
Then, \( f_x\dot{x}+ f_y\dot{y}+ f_z\dot{z}=
0
\) implies \(\alpha_{m,n,k}=0\) for all \(m,n,k\in \NZ,\) \ie \(f=0.\)
This completes the proof.
\epr

Now we briefly present a general framework on how to obtain the simplest (orbital and parametric) normal form of a vector field, for more detailed information, terminology and background material see \cite{GazorYuSpec}. Let \(\E=\sum^\infty_{i=0} \E_i\) be a graded Lie algebra, \ie \([v_i,v_j]\in \E_{i+j}\) for any \(v_i\in \E_i\) and \(v_j\in\E_j\). The sub-indices \(i\) for spaces (not element indices, of course) in this paper refer to the homogenous spaces of grade \(i\) according to the associated grading. Recall that a style of a normal form is a rule to choose a unique complement space to a subspace of \(\E.\) Assume a normal form style is given and \(A= \sum^\infty_{i=1}A_i\) is the graded vector space of all generators associated with the permissible transformations (in \cite{GazorYuSpec} this is denoted by \(A^0\)). The space \(A\) acts linearly on \(\E.\) This action is denoted by \(*\) and it preserves the grading structures, \ie \(A_i*\E_j\subseteq \E_{i+j}\) for any \(i, j.\)
\begin{rem}
The action \(*\) is the initially linear map of the associated near-identity permissible transformations generated by elements of \(A\), see \cite{GazorYuSpec} for more detailed discussion.
\end{rem}
The space \(A\) and its associated action \(*\) are different for the normal form, orbital normal form and parametric normal form of vector fields. Therefore, we will explicitly define these accordingly in the following three sections. Let \(v=\sum^\infty_{i=0}v_i\) and inductively define
\be
d^{n,1}: A_n\rightarrow\E_n, \quad \hbox{ by } d^{n,1}(Y_n):= Y_n*v_0 \hbox{ for } Y_n\in A_n,
\ee and
\be
d^{n,k}: A_n\times \ker(d^{n-1,k-1})\rightarrow\E_n \quad (\hbox{for any } k\leq n),
\ee given by \(d^{n,k}(Y_{n}, Y_{n-1}, \ldots, Y_{n-k+1}):= \sum^{k-1}_{i=0} Y_{n-i}*v_i,\) where
\bes (Y_{n-1}, Y_{n-2}, \ldots, Y_{n-k+1})\in \ker(d^{n-1,k-1}).\ees
For any \(k>n,\) let \(d^{n,k}:=d^{n,n}.\)
Note that \((Y_{n-1}, Y_{n-2}, \ldots, Y_{n-k+1})\in \ker(d^{n-1,k-1})\) if and only if
\be\sum^{j}_{i=0} Y_{n+i-j-1}*v_{j-i}=0\ee for any \(j=0, 1, \ldots,k, \) see also \cite{wlhj}. Let
\(\RC^{n,k}:= \IM (d^{n,k}).\) Then, there exist the complement subspaces \(\C^{n,k}\) such that
\be
\RC^{n,k}\oplus \C^{n,k}=\E_n
\ee follows the normal form style. We call \(d^{n,k}\) the \(k\)-th level state map.
If \(v_i\in \C^{n,k}\) for any natural number \(n,\) we say that \(v=\sum v_i\) is in the \(k\)-th level extended partial normal form (see \cite{Mord04} for the original definition and ideas) and when \(v_i\in \C^{i,i}\) for any natural number \(i,\) \(v\) is said to be in the infinite level (orbital or parametric) normal form. Notice that the domain of the map \(d^{n,k}\) has \(k\)-components. However, the first component does not play any role in our normal form computation, since \(\Theta^0_0\) is the zero grade of the system and we have \([\E, \Theta^0_0]=0\) in Section \ref{secNF}, and in Sections \ref{secONF} and \ref{secPNF} the linear part \(\Theta^0_0\) is first removed before normal form computation. Thus, we omit the first component and we only have \((k-1)\)-components in the domain of the map \(d^{n,k},\) \ie \(d^{n,k}(Y_{n-1}, \ldots, Y_{n-k+1})= \sum^{k-1}_{i=1} Y_{n-i}*v_i.\)

\begin{thm}\label{LNF}
There always exist a sequence of permissible transformations transforming a vector field into its \(k\)-th (or infinite) level extended partial normal form.
\end{thm}
\bpr
The proof follows from \cite[Theorem 4.4]{GazorYuSpec}.
\epr

In this paper, we also discuss the symmetry group of the infinite level normal forms. The symmetry group of a vector field is a subgroup of the permissible transformations such that the infinite level normal form \(v^{(\infty)}\) is invariant under the symmetry group acting on \(v^{(\infty)}\). This symmetry group is generated by the infinite level co-normal form space. The infinite level co-normal form space is a subspace, say \(B,\) of the space of permissible transformation generators \(A\) in which \(B*v^{(\infty)}=0\) and can be computed through the normal form computation.

\section{The simplest normal forms }\label{secNF}

In this section we obtain the simplest normal form of vector fields given by Equation (\ref{Eq1}) with nonzero quadratic terms, \ie
\be\label{E0}
v^{(1)}:= \Theta^0_0+ \beta^1_1 E^1_1 + \cdots\in \Theta^0_0+\E,
\ee where \(\beta^1_1\neq0.\) Define the grading function \(\delta \) for the space \(\E\) by
\bas
\delta(E^l_k)&:=& 2k-l \quad \hbox{ and } \quad \delta(\Theta^0_0)=0.
\eas Then, \(\delta(E^l_k)\geq k.\) For simplicity of notation \(\delta\) or \(\delta_\alpha\) is used for all gradings in this paper.

\begin{lem}\label{seclev}
There exists a formal near identity transformation that sends the vector field \(v\) given by Equation (\ref{E0}) into the second level extended partial normal form
\ba\label{secNFE}
v^{(2)}=\Theta^0_0+\beta^1_1 E^1_1+\beta^0_{1}E^0_1+\beta^2_2E^2_2+\beta^1_2E^1_2+\sum^\infty_{k=2}\beta^0_kE^0_k,
\ea where \(\beta^1_1\) is not changed.

\end{lem}
\bpr Define \(A:= \sum^\infty_{i=1}\E_i.\) Since \([E^l_k, E^{1}_{1}]= -(k-1)E^{l+1}_{k+1},\) a complement space to the image of \(d^{n,2}\) is \(\Span \{E^0_{2k}\}\) for any \(n=2k>3,\) \(\Span \{E^1_{2}\}\) for \(n=3,\) and \(\Span \{E^2_{2}, E^0_{1}\}\) for \(n=2.\)
\epr

The reader should note that in this paper the notation \(\beta^l_k\) for coefficients are not changed as they are updated in the normal form computation.

\begin{thm}\label{SNF} There exists a natural number \(r\) and a
formal near identity transformation sending \(v^{(2)}\) given by equation (\ref{secNFE}) into one of the following \((2r+2)\)-th level extended partial normal forms:
\be\label{2n2f}
v^{(2r+2)}= \Theta^0_0+\beta^1_1 E^1_1+\beta^0_{1}E^0_1+\beta^2_2E^2_2+\beta^1_2E^1_2+ \beta^0_rE^0_r +\sum^\infty_{k=r+2}\beta^0_kE^0_k \quad(\hbox{for } r\geq 2),
\ee
or
\be\label{2n2s}
v^{(2r+2)}= \Theta^0_0+\beta^1_1 E^1_1+\beta^0_{1}E^0_1+\beta^2_2E^2_2+\beta^1_2E^1_2+ \beta^0_1(\beta^1_2-\beta^0_1\beta^2_2)E^0_2+ \beta^0_rE^0_r +\sum^\infty_{k=r+2}\beta^0_kE^0_k,
\ee where \(\beta^0_1(\beta^1_2-\beta^0_1\beta^2_2)\neq0\)  and \(r>2.\)
Furthermore, \(v^{(2r+2)}=v^{(\infty)},\) and for \(\beta^0_r\neq 0\) there is no non-trivial symmetry in the symmetry group for \(v^{(\infty)}\).
\end{thm}
\bpr
Consider the vector field \(v^{(2)}\) obtained in Lemma \ref{seclev}. Denote
\bes\mathcal{E}_{0,k+1}:=\Big(0,\ldots,0,(\beta^1_2-\beta^0_1\beta^2_2)E^0_2,\beta^2_2E^1_2,E^0_1,0\Big)\ees
and \bes\mathcal{E}_{1,k+2}:=\Big(0, \ldots, 0, \big(\beta^0_2-\beta^0_1(\beta^1_2-\beta^0_1\beta^2_2)\big)E^0_2, (\beta^1_2-\beta^0_1\beta^2_2)E^1_2, \beta^2_2E^2_2, E^1_1\Big),\ees
where the number of zeros is \((k-3)\)-times for any \(k\geq 3.\)
Then, we have
\bas
\ker(d^{5,5})= \Span \{\mathcal{E}_{0,5}\}.
\eas

In the case that \(\beta^0_2-\beta^0_1(\beta^1_2-\beta^0_1\beta^2_2)\neq 0,\) we have
\bes d^{6,6}(\mathcal{E}_{0,6})=\big(\beta^0_2-\beta^0_1(\beta^1_2-\beta^0_1\beta^2_2)\big) E^0_3\in \IM d^{6,6}
\hbox{ and } d^{6,6}(\mathcal{E}_{1,6})=\beta^0_1\big(\beta^0_1(\beta^1_2-\beta^0_1\beta^2_2)-\beta^0_2\big) E^0_3.\ees
Then,  \(\beta^0_1\mathcal{E}_{0,6}+\mathcal{E}_{1,6}
\in \ker d^{6,6}.\) Let \(r:=2,\) and hence the \(6\)-th level extended partial normal form is given by Equation (\ref{2n2f}).

For \(\beta^0_2=\beta^0_1(\beta^1_2-\beta^0_1\beta^2_2),\) let \(r:= \min \{k\,|\; \beta^0_k\neq 0 \hbox{ for } k>2\}.\)
When \(\beta^0_k= 0\) for all \(k>1,\) the choice \(r:=2\) trivially completes the proof. Thus, we can assume that \(\beta^0_r\neq 0\) and \(r\) is finite. Therefore, \(\ker d^{k,k}=\Span \{\mathcal{E}_{0,k}, \mathcal{E}_{1,k}\}\) for any \(6\leq k\leq 2r+1,\)
\bas
d^{2r+2,2r+2}(\mathcal{E}_{0,2r+2})=(r-1)\beta^0_{r}E^0_{r+1},
& \hbox{ and } &
d^{2r+2,2r+2}(\mathcal{E}_{1,2r+2})={-(r-1)\beta^0_{r}\beta^0_1E^0_{r+1}}.
\eas Thereby, \bes E^0_{r+1}\in \IM d^{2r+2,2r+2} \hbox{ and } \ker d^{2r+2,2r+2}=\Span\{\mathcal{E}_{2r+2}\},\ees where \(\mathcal{E}_{2r+2}:=\beta^0_1\mathcal{E}_{0,2r+2}+ \mathcal{E}_{1,2r+2}.\)
Hence, \(\ker d^{k,k}=\Span\{\mathcal{E}_{k}\}\) for \(k\geq2r+2\) is one dimensional. Identify \(\mathcal{E}_{k}\) with an element in \(\E,\) then \(\mathcal{E}_{k}\) converges to the infinite level normal form vector field in the filtration topology when \(k\) approaches infinity, see \cite{GazorYu,GazorYuFormal}. This proves that the simplest normal form computation does not produce any non-trivial symmetry for the normal form system. The proof is complete by Theorem \ref{LNF}.
\epr

\section{The simplest orbital normal form}\label{secONF}

In this section, we obtain the simplest orbital normal form of \(v\) given by Equation (\ref{E0}).
Let \(Z^m_n:= x^my^{n-m}\) for \(m\leq n, m, n\in \N,\) and \(Z^0_0:=1.\)
Then,
\bas
\RM:= \Span \Big\{\sum_{n\geq m\geq 0} a^m_nZ^m_n| a^m_n\in \f, n\neq 0\Big\}
\eas
is the space associated with time rescaling generators, \ie \(\f[[x, y]]= \RM\oplus \f.\) The time rescaling computation is governed by a module structure, see \eg \cite{GazorYuSpec}.
Thus, the structure constants for the ring \(\RM\) and the \(\RM\)-module \(\E\) are given by
\ba
Z^m_nZ^l_k:=Z^{m+l}_{n+k} &\hbox{ and }&  Z^m_n E^l_k:= E^{m+l}_{n+k}.
\ea Furthermore, let
\be
r:= \min \{k|\beta^0_{k}\neq0\}.
\ee A grading function is defined by \(\delta(Z^l_k):= 2k-l\) and thus, \(\RM= \sum^\infty_{i=1}\RM_i.\) Define \(A_i:=\E_i\times \RM_i\) for any \(i\geq 1,\) \(A:= \sum^\infty_{i=1}A_i.\) Here \(A\) is a graded vector space according to \cite{GazorYuFormal,GazorYuSpec} and denotes the set of all formal sums \(\sum^\infty_{i=1}(v_i, T_i),\) where \(v_i\in \E_i\) and \(T_i\in \RM_i\). The action of \(A\) on \(\E\) is given by \((S,T)*v:=[S,v]+Tv\) for any \((S, T)\in A\) and \(v\in \E.\)

\begin{thm}\label{inftThmONF}
There exists a natural number \(r\) such that the vector field \(v\in \Theta^0_0+\E\) with nonzero quadratic part can be transformed into the simplest orbital normal form
\bas\label{inftONF}
v^{(\infty)}&=&\Theta^0_0+ E^1_1+\beta^0_{1}E^0_1+\beta^0_r E^0_r, \ \hbox{ for any } 
r\geq2.
\eas Furthermore, the symmetry group of the vector field \(v^{(\infty)}\) is not finitely generated. When \(\beta^0_r\neq 0\) the symmetry group has no non-trivial symmetry from the group of permissible changes of state variables.
\end{thm}
\bpr We first remove the linear part \(\Theta^0_0\) from the system with a linear change of state variable \([x,y,z]:=exp\big((Z\ddY-Y\ddZ)t\big)[X,Y,Z],\) that sends the new variables \([X,Y,Z]\) into the old variables, see also \cite[Proof of Lemma 5.3.6]{MurdBook}. Therefore, we may assume \(v:=E^1_1+\ldots\in \E.\) Since \bes d^{2,2}(0,Z^1_1)= -\beta^2_2E^2_{2}\in \IM d^{2,2}\hbox{ and } d^{3,2}(0,Z^0_1) =-(\beta^1_2-\beta^2_2\beta^0_1) E^1_{2}\in \IM d^{3,2},\ees both terms of \(E^2_{2}\) and \(E^1_{2}\) are eliminated from the second level extended partial orbital normal form, \ie \(v^{(2)}= E^1_1+\sum^\infty_{k=1}\beta^0_k E^0_k.\) Now let \(r:=\min \{k| \beta^0_k\neq0, k>1\}.\) (In case \(r\) is not finite, an arbitrary choice for \(r>1\) trivially completes the proof.)
Denote \bes \ZS^0_{m,k+2}:= \big(0,\ldots,0, Z^0_m, \frac{-1}{(1-m)}E^0_m\big)\ees where the number of zeros is
\(2k\)-times for \(k\geq 1\). Then,
\be\label{2rl}
\ZS^0_{m,3}= \Big(0,0, Z^0_m, \frac{-1}{(1-m)}E^0_m\Big)\in \ker d^{2m+2,3}\hbox{ and }
d^{2m+2r,2r+1}(\ZS^0_{m,r})=\frac{r-1}{m-1}\beta^0_r E^0_{m+r},
\ee
for any \(m\geq 1.\)
By Theorem \ref{LNF} there exist transformations that send the vector field \(v\) into the infinite level orbital normal form \(u:=E^1_1+\beta^0_{1}E^0_1+\beta^0_r E^0_r\in \E\) for some \(r\geq2, \beta^0_1, \beta^1_r\in \f.\)
For any \(Z^l_k\) (\(l\geq 1, k\geq 2)\) there exists a nonzero state solution \(S\in \mathscr{S}\) and a time solution \(T\in \ta\) such that \((T, S)*u=-Z^l_ku. \) Therefore, the symmetry group for \(u\) is not finitely generated. Now the linear change of state variable \([x,y,z]:=exp\big((Y\ddZ-Z\ddY)t\big)[X,Y,Z]\) puts the linear part \(\Theta^0_0\) back into the normal form system, \ie \(v^{(\infty)}\) in Equation (\ref{inftONF}) is obtained. The composition of this map with the symmetry group associated with \(u\) results in a symmetry group for \(v^{(\infty)}.\) The rest of the proof follows Theorem \ref{SNF}.
\epr

\section{Parametric normal form}\label{secPNF}

This section is concerned about parametric normal form of a non-degenerate multi-parameter nonconservative perturbation \(w(\mu)\) (will be defined bellow) of the vector field \(v\) given by Equation (\ref{Eq1}), see also \cite{GaoPNF,GazorYuSpec,Mord08}. Define
\be\label{LieP}
\mathscr{S}:= \Span\Big\{\sum\beta^l_{k,\mathbf{m}} E^l_k\mu^\mathbf{m}| \beta^l_{k,\mathbf{m}}\in \f\Big\}, \ee where
the summation is over \(\mathbf{m}\in \NZ^q, \) \({l+k\geq 0},\) \(k\geq l,\) \(\mu:= (\mu_1, \mu_2, \ldots, \mu_q), |\mathbf{m}|=\sum^q_{j=1}m_j\geq 0, \beta^0_{0,\0}=0,\) and \(\mu^\mathbf{m}:= \mu_1^{m_1}\ldots\mu_q^{m_q}.\) We similarly define the parametric time space \(\ta\) such that \(\f[[x,y,\mu]]=\ta\oplus \f,\) and the space \(\PM\) as the space of all formal power series in terms of parameters \(\mu\in \f^q\) without constant term,
\ie \(\f[[\mu]]^q=\PM\oplus \f^q.\) The algebraic structures are naturally extended from nonparametric to the parametric cases,
see \eg \cite{GazorYuFormal,GazorYuSpec}. For \(\mathbf{m}=(m_1, \ldots, m_q)\in \NZ^q\) and \(\alpha=(\alpha_1, \ldots, \alpha_q)\in \N^q,\)
let \(\alpha\cdot\m:=\sum^q_{j=1}\alpha_jm_j.\) Then, define the grading function by
\bes
\delta_\alpha(E^l_k\mu^\mathbf{m}):= 2k-l+\alpha\cdot\mathbf{m}, \delta_\alpha(Z^l_k\mu^\mathbf{m}):= 2k-l+\alpha\cdot \m,\hbox{ and }
\delta_\alpha(\mu^\mathbf{m}):=\alpha\cdot \m.\ees
Define \bes A_i:= \mathscr{S}_i\times\PM_i\times \ta_i,\; A:=\sum^\infty_{i=1} A_i, \hbox{ and } (S, P, T)*v:= D_\mu (v)P+ Tv+[S, v]\ees for any \((S, P, T)\in A.\) The perturbation of the system is considered to be within the space of \(\mathscr{S},\) \ie \(w(\mu)\in \Theta^0_0+\mathscr{S}\). We call any vector field
\be\label{pv} w(\mu)\in \Theta^0_0+ \mathscr{S}\ee
a multi-parameter nonconservative perturbation of the vector field \(v\) given by Equation (\ref{Eq1}), when \(w(0)=v.\) Similar to the previous section, we remove the linear part \(\Theta^0_0\) from the vector field with a linear change of state variable. Since \(w_0=0,\) the first component of the map \(d^{n,k}\) is omitted, \eg the map
\be
d^{n,2}: A_{n-1}\rightarrow\E_n,
\ee is defined by \(d^{n,2}(Y_{n-1})= Y_{n-1}*v_1 \hbox{ for any } Y_n\in A_n.\) Since \(\beta^1_{1,\mathbf{0}}\neq 0,\) through a time rescaling we may assume that \(\beta^1_{1,\mathbf{0}} =1,\) \ie \(w=\sum^\infty_{i=1}w_i,\) \(w_1=E^1_1\) and \(w_i= \sum \beta^l_{k,\mathbf{m}} E^l_k\mu^\mathbf{m}\) with the summation over \(l, k, \mathbf{m}\in \NZ^q,\) \(i=2k-l+\alpha\cdot\mathbf{m}\) for any \(i>1.\)

\begin{lem} Any multi-parameter nonconservative parametric perturbation of a vector field given by equation (\ref{Eq1}) can be transformed to the second level extended partial  parametric normal form of \(w^{(2)}\) given by
\be\label{sec}
w^{(2)}= E^1_1+ \sum_{k\geq 0, |\mathbf{m}|\geq 0} \beta^0_{k, \mathbf{m}}E^0_k \mu^{\mathbf{m}} \quad (\beta^0_{0,\0}=0).
\ee
\end{lem}
\bpr Let \(\alpha:=(1,\ldots,1).\) Since \bes D_\mu v_1=0,\, [E^l_k,E^1_1]= (1-k)E^{l+1}_{k+1},\,
Z^0_0E^1_1= E^1_1, \, Z^0_1E^1_1= E^1_2, \,\hbox{ and } Z^1_1E^1_1= E^2_2,\ees for any \(N\geq 2\) we have
\be\IM d^{N,2}= \Span \{E^l_k\mu^{\mathbf{m}}\,| \hbox{ for any } l\geq 1, k\geq l, \hbox{ where } N= 2k-l+|\mathbf{m}|\}.\ee
Thereby, \(\NC^{N,2}= \Span\{E^0_k\mu^{\mathbf{m}}\,|\; N= 2k+|\mathbf{m}|\}\)  for any \(N\geq 2.\)
\epr
Consider the second level extended partial parametric normal form \(w^{(2)}\) given by Equation (\ref{sec}). Let \(r:=\min \{k\,|\; \beta^0_{k,\mathbf{0}}\neq0, k>1\}.\) For the rest of this paper, we assume that there exists \(\beta^0_{k,\mathbf{0}}\neq0\) for a \(k>1,\) \ie \(r\) is finite. Define a \((r+1)\times q\) matrix
\be
A:= (a_{ij}), \hbox{ where } a_{ij}:= \beta^0_{i-1,e_j} \hbox{ and } e_j \hbox{ denotes the } j\hbox{-th standard basis of } \f^q.
\ee

\noindent Then, we call a parametric vector field \(w(\mu)\in \Theta^0_0+\mathscr{S}\) a non-degenerate multi-parameter nonconservative perturbation of the vector field \(v\) if \(v=w(\mathbf{0}),\) and \(\rank(A)= r+1\) when \(w(\mu)\) is transformed into \(w^{(2)}\) given by Equation (\ref{sec}). When this condition is satisfied, through an invertible linear change of parameters we can transform the vector field \(w^{(2)}\) into
\be\label{secT}
\tilde{w}^{(2)}= E^1_1+ \sum^r_{k=1}\beta^0_{k,\mathbf{0}}E^0_k+ \sum^r_{k=0}E^0_k \mu_{k+1}+ \sum_{r\geq k\geq 0, |\mathbf{m}|\geq 2} \beta^0_{k, \mathbf{m}}E^0_k\mu^{\mathbf{m}} + \sum_{k> r, |\mathbf{m}|\geq 0} \beta^0_{k, \mathbf{m}}E^0_k \mu^{\mathbf{m}},
\ee
\ie \(a_{ij}=0\) if \(i\neq j,\) and \(a_{ii}=1.\) Now we update the grading weight vector \(\alpha\) by \(\alpha_i:=2(r-i+1).\) Therefore, \(\delta_\alpha(E^0_k \mu_{k+1})= 2r\) for any \(0\leq k\leq r.\)

\begin{lem} \label{SPNF}Let \(w(\mu)\) be a non-degenerate multi-parameter nonconservative perturbation of the vector field \(v\) governed by Equation (\ref{Eq1}). Then, there exist a natural number \(r\) and a sequence of permissible transformations such that they transform \(w\) into its \((2r+1)\)-th level extended partial parametric normal form
\be
w^{(2r+1)}= E^1_1+ \beta_{1}E^0_1+ \beta_{r}E^0_r+ \sum^r_{k=0} E^0_k \mu_{k+1} \quad (\beta^0_{r,\mathbf{0}}\neq 0).
\ee Furthermore, \(w^{(\infty)}=w^{(2r+1)}.\) The symmetry group of \(w^{(\infty)}\) is not finitely generated, yet it has no non-trivial symmetry from the group of permissible changes of state variables.
\end{lem}

\bpr Without loss of generality we assume that  \(w\) is transformed into \(\tilde{w}^{(2)}\) given by Equation \ref{secT} and \(\beta^0_{r,0}\neq0\).
By Equation (\ref{2rl}), we have
\bes E^0_{m+r}\mu^\mathbf{m}\in \IM d^{2m+2r+\alpha\cdot\mathbf{m}, 2r+1} \ \hbox{ for any } m\geq 1.\ees Furthermore, the condition \(\rank(A)=r+1\) implies that \bes E^0_k\mu^\mathbf{m}\in \IM d^{2k+\alpha\cdot\mathbf{m}, 2r+1} \ \hbox{ for any } |\mathbf{m}|\geq 2 \hbox{ and } 0\leq k\leq r.\ees Theorem \ref{LNF} implies that there exist transformations such that \(w^{(2)}\) can be transformed into \(w^{(2r+1)}\) as desired. Since \(w^{(2r+1)}\) does not have any term of grade higher than \(2r,\) we have \(w^{(\infty)}=w^{(2r+1)}.\) The uniqueness of \(w^{(\infty)}\) together with an argument similar to the proof of Theorem \ref{inftThmONF} complete the proof.
\epr
A linear change of state variables creates the linear part \(\Theta^0_0\) back into the infinite level normal form system. This together with Lemma \ref{SPNF} proves the following theorem.
\begin{thm}\label{PNF} Let the hypothesis of Lemma \ref{SPNF} hold. Then, there is a sequence of permissible transformation transforming the system generated by \(w(\mu)\) into the infinite level normal form system
\ba\label{FinalPNF}
{\dot{x}}&=&x g(x,\rho,\mu), \quad {\dot{\rho}}=\frac{1}{2}\rho g(x,\rho,\mu), \quad \dot{\theta}=1,
\ea where \(g(x,\rho, \mu)=x+\beta_1\rho^{2}+\sum^{r}_{k=0}\rho^{2k}\mu_{k+1}+\beta_r \rho^{2r},\) \(\beta_r\neq0\) and \(r\geq2.\)
\end{thm}

The following theorem sheds light to the dynamics of the infinite level parametric normal form system.
\begin{thm}\label{bifur} Consider the infinite level parametric normal form system in Theorem \ref{PNF}. Then, the surface \(\rho^2=c|x|,\) for any nonnegative constant \(c\), is an invariant manifold of solutions for the system. Furthermore, the surface \(x=-\beta_1\rho^{2}-\sum^{r}_{k=0}\rho^{2k}\mu_{k+1}-\beta_r \rho^{2r}\) is an invariant manifold of limit cycles for the parametric system; except that it includes the origin as a fixed point. This implies that any point outside the later invariant surface, approaches a limit cycle or the origin through the surface \(\rho^2=c|x|\) for an appropriate constant \(c.\)
\end{thm}
\bpr By Equation (\ref{Eul}) we have \(\frac{1}{2}\rho\dot{x}= x\dot{\rho}\) and thus, \(\rho^2=c|x|\) is an invariant set for the system.
The rest of the proof follows from Theorem \ref{PNF}.
\epr

Now we present an example and find its parametric normal form.
\begin{exm} Consider the system governed by
\bas
\dot{x}&=& 2x\mu_1+2x^2(a+\mu_{1}+\mu_{2})+2x(y^2+z^2)(b+\mu_{2})+2x^3(c+{\mu_{1}}^2+{\mu_2}^2)+ 2x(y^2+z^2)^2\mu_3
\\&&+ 2ex^4+2x^2(y^2+z^2)(d+\mu_{1}^2),
\\
\dot{y}&=&z+y\mu_1+xy(a+\mu_{1}+\mu_{2})+y(y^2+z^2)(b+\mu_{2})+x^2y(c+{\mu_2}^2)+ y(y^2+z^2)^2\mu_3
\\&&+ex^3y+xy(y^2+z^2)(d+\mu_{1}^2),
\\
\dot{z}&=&-y+z\mu_1+xz(a+\mu_{1}+\mu_{2})+z(y^2+z^2)(b+\mu_{2})+x^2z(c+{\mu_2}^2)+ z(y^2+z^2)^2\mu_3
\\&&+ex^3z+xz(y^2+z^2)(d+\mu_{1}^2).
\eas where \(ab(da-cb)\neq 0.\) Then (using a Maple program) we obtain its parametric normal form given by
\bas
\dot{x}&=& 2x g(x, y^2+z^2, \mu),
\\
\dot{y}&=&z+y g(x, y^2+z^2, \mu),
\\
\dot{z}&=&-y+z g(x, y^2+z^2, \mu),
\eas where \(g(x, R, \mu):=\mu_1+(\beta_1+\mu_2)R+(\beta_2+\mu_3) R^{2},\) \(\beta_1=\frac{b}{a}\) and \(\beta_2=\frac{b(da-cb)}{a^3}.\)

\end{exm}


\end{document}